\documentclass[12pt]{amsart}
\newtheorem{theorem}{Theorem}[section]
\newtheorem{lemma}[theorem]{Lemma}

\theoremstyle{definition}

\theoremstyle{remark}
\newtheorem{remark}[theorem]{Remark}

\theoremstyle{conj}

\numberwithin{equation}{section}

\begin{document}
\title{A note on the Hitchin-Thorpe inequality and Ricci flow on 4-manifolds  }

\author[Y. Zhang]{Yuguang Zhang*}
\thanks{*Supported by  NSFC-10901111, and KM-210100028003. }
\address{Department of Mathematics, Capital Normal University,
Beijing, P.R.China  }
 \email{yuguangzhang76@yahoo.com}
\author[Z. Zhang]{Zhenlei Zhang**}
\thanks{**   Supported by  NSFC-09221010056}
\address{Department of Mathematics, Capital Normal University,
Beijing, P.R.China }
 \email{zhleigo@yahoo.com.cn}

\begin{abstract}
In this short  paper, we prove a Hitchin-Thorpe type  inequality for
closed 4-manifolds with non-positive   Yamabe invariant, and
admitting long time solutions  of the normalized Ricci flow equation
with bounded scalar curvature.
\end{abstract}

\maketitle

\section{Introduction}
A Riemannian metric $g$ on a smooth manifold $M$  is called an
Einstein metric,  if   $$ Ric_{g}=cg,$$ where $ Ric_{g}$ is the
Ricci tensor, and $c$ is a constant.  If a closed  oriented
4-manifold $M$ admits an Einstein metric $g$, there is an
inequality, the Hitchin-Thorpe inequality,  for the Euler number
$\chi(M)$  and the signature $\tau (M)$  of $M$:
  \begin{equation}\label{*0} 2\chi(M)-3|\tau (M)|\geq 0,  \end{equation}
 (cf.  \cite{Hi}, or Theorem 6.35 in \cite{Be}). This inequality
 serves as a topological obstruction for the existence of Einstein
 metrics on 4-manifolds, i.e. if (\ref{*0}) is not satisfied, then
 $M$ wouldn't admit any Einstein
 metric.

  The   Ricci flow was introduced by Hamilton in \cite{H}  to find Einstein metrics on a
  given manifold, which   is the following  evolution equation for a  smooth family
of metrics $g(t),t\in[0,T),$ \begin{equation}\label{00*}
\frac{\partial}{\partial t}g(t)=-2Ric_{t}.
\end{equation} The normalized Ricci flow is
  \begin{equation}\label{*}
\frac{\partial}{\partial t}g(t)=-2Ric_{t}+\frac{2r(t)}{n}g(t)
\end{equation}
where $r(t)=\frac{\int_M R_{t}dv_{g(t)}}{{\rm Vol}_{g(t)}(M)}$
denotes the average scalar curvature, and $R_{t}$ denotes the scalar
curvature of $g(t)$.  The normalized Ricci flow  is just a
transformation of (\ref{00*}) by rescaling the space and time such
that the volume preserves to be a constant along the flow. If
(\ref{*})  admits a long time solution $g(t),t\in[0,\infty)$, and
$g(t)$ converges to a Riemannian metric $g_{\infty}$ on $M$ in some
suitable sense, when $t\rightarrow\infty$, then $g_{\infty}$ is an
Einstein metric or a Ricci-soliton.  However, (\ref{*}) may not
admit  any  long time solution, and, even  there is such a solution,
$g(t)$ may not converge to a Riemannian metric  on $M$. It is
expected that the inequality  (\ref{*0}) is also  a topological
obstruction for the existence of long time solutions of (\ref{*}) on
4-manifolds, at least for 4-manifolds with non-positive Yamabe
invariant. In \cite{FZZ1}, an analog inequality of (\ref{*0}) was
obtained for long times solutions of (\ref{*}) under some
hypothesis. In this paper, we will continue to study the
Hitchin-Thorpe type inequality for 4-manifolds admitting long time
solutions of (\ref{*}).

For  a closed   Riemannian $n$-manifold  $(M,g)$ and a function
$f\in C^{\infty}(M)$, set
$$\mathcal{F}(g,f)=\int_{M}(R_{g}+|\nabla f|^{2})e^{-f}dv_{g}. $$   The Perelman's
$\lambda$-functional is defined by
$$\lambda_{M}(g)=\inf _f \{\mathcal{F}(g,f)| \int_{M}e^{-f}
 dv_{g}=1\},$$
  which  is the lowest eigenvalue of the
 operator $-4\triangle+R_{g}$.    Perelman \cite{P1} has
established the monotonicity property of $\overline{\lambda}
_M(g(t))=\lambda_{M}(g(t))\text{Vol}_{g(t)}(M)^{\frac{2}{n}}$ along
the Ricci flow. A  diffeomorphism invariant
  $\overline{\lambda}_{M} $ of $M$ is defined  (cf. \cite{P2}  \cite{KL})
by $$\overline{\lambda}_{M}=
  \sup\limits_{g\in \mathcal{M}} \overline{\lambda}_{M}(g),
  $$ where $\mathcal{M}$ is the set of Riemannian
metrics on $M$.   By \cite{AIL} and \cite{An1},   $\overline{\lambda
}_M$ is equal to the classical   Yamabe invariant whenever
$\overline{\lambda}_{M}\leq 0$ or the Yamabe invariant is
nonpositive.

 In \cite{FZZ1}, it was  proved that,
if $g(t), t\in[0,\infty)$, is a solution to (\ref{*}) with
$|R(g(t))|<C$ for a constant $C$ independent of $t$ on a   closed
oriented 4-manifold $M$ with $\overline{\lambda}_{M}< 0$,  then
$$2\chi(M)-3|\tau(M)|\geq
\frac{1}{96\pi^{2}}\overline{\lambda}_{M}^{2}>0.$$
  Based on this inequality and the  Seiberg-Witten theory, Ishida \cite{I} recently
showed that the existence of long time non-singular solution really
depends on the smooth structure of the underlying manifold.
 The goal of the present paper is to  relax the assumption $\overline{\lambda}_{M}< 0$ to
 $\overline{\lambda}_{M}\leq
 0$.

 \vskip 5mm

\begin{theorem}\label{1.1} Let $M$ be a closed oriented  $4$-manifold with $\overline{\lambda}_{M}\leq 0 $.  If  $M$  admits a long time  solution
 $g(t),t\in[0,\infty),$ of (\ref{*})  with scalar curvature  $|R_{t}|<C $ for a constant $C$
 independent of $t$,  then
$$2\chi(M)-3|\tau (M)|\geq \frac{1}{96\pi^{2}}\overline{\lambda}_{M}^{2}.
$$
\end{theorem}

 \vskip 5mm

 The hypothesis of $\overline{\lambda}_{M}\leq 0 $
    holds for many cases, for example,    complex projective
    surfaces
    with non-negative  Kodaira  dimension by  \cite{L}, which include  K3-surfaces,  projective  surfaces of general type
    and some  surfaces of
    elliptic type etc.
  Furthermore, Theorem 2 in  \cite{Ko}
 shows  that, if $M$ is a closed oriented  $4$-manifold with a monopole class $c_{1}(\mathfrak{c})$ that is not a
  torsion class and satisfies $c_{1}^{2}(\mathfrak{c})\geq 0$, then $$ \overline{\lambda}_{M}\leq
  -\sqrt{32\pi^{2}c_{1}^{2}(\mathfrak{c})}\leq 0$$ (see also \cite{FZZ2} for the  case of $c_{1}^{2}(\mathfrak{c})> 0$).
   In \cite{FZZ1},    a
 Miyaoka-Yau type
  inequality was proved for 4-manifolds with a   monopole class
  $c_{1}(\mathfrak{c})$ such that   $c_{1}^{2}(\mathfrak{c})> 0$,  and admitting long time  solutions  of (\ref{*})  with bounded  scalar
  curvature.

The assumption of bounded scalar curvature is
       a technique assumption, and we   hope that it can be removed in the
       future study. However,  it can be verified for some cases.
       In \cite{Cao}, it was shown that a K3-surface $M$ admits a long
       time solution   $g(t)$, $t\in [0,\infty) $,  of (\ref{*}), and $g(t)$ converges
       to a
       Ricci-flat K\"ahler-Einstein  metric  on $M$. Thus the scalar
       curvature of $g(t)$ is uniformly bounded. Since $\overline{\lambda}_{M}=0$ (cf.
       \cite{L}), Theorem \ref{1.1} can apply to this case.
  If $M$ is a  complex minimal projective  surface of  elliptic  type
    with   Kodaira  dimension 1 and no singular fibers, then
    $\overline{\lambda}_{M}=0$ by \cite{L}, and $M$  admits
        long time  solutions $\tilde{g}(\tilde{t})$, $\tilde{t}\in [0,\infty) $,  of the following K\"ahler-Ricci flow equation
   \begin{equation}\label{10*}\frac{\partial}{\partial \tilde{t}}\tilde{g}(\tilde{t})=
-Ric_{\tilde{t}} -\tilde{g}(\tilde{t})  \end{equation}  by \cite{Ts}
and  \cite{TZ}. By Corollary 1.1 in  \cite{ST},  the scalar
curvature of $\tilde{g}(\tilde{t})$ is uniformly  bounded along
(\ref{10*}). A straightforward calculation would  show that  a
transformation of
 $\tilde{g}(\tilde{t})$, $\tilde{t}\in [0,\infty) $,   by rescaling the metric  and
 time gives a long time solution of (\ref{*}) with bounded scalar
 curvature.  Thus the assumption  of  Theorem \ref{1.1} is
 satisfied. 
  If $M$ is a complex minimal projective  surface of
    general type, which satisfies  $\overline{\lambda}_{M}<0$ by
    \cite{L}, then
    $M$ admits   long time solutions  of (\ref{*}) with bounded scalar
 curvature by \cite{Cao},  \cite{Ts},  \cite{TZ}, \cite{ZZ} and the same arguments as above.    In
\cite{ZY}, an alternative proof of the  Miyaoka-Yau
  inequality for  minimal
projective manifolds of general type was obtained by  using  the
result in \cite{ZZ} and a similar Ricci flow argument.

  An analog  Hitchin-Thorpe type inequality for non-compact 4-manifolds admitting
  non-singular solutions of (\ref{*}) is  obtained in \cite{FZZ3}.   We assume ${\rm Vol}_{g(t)}(M) \equiv 1$ in this paper for
convenience. We shall prove   Theorem \ref{1.1}  in the  next
section.

\vskip 10mm

\section{Proof of Theorem \ref{1.1}}
The goal of this section is to prove  Theorem \ref{1.1}, which
essentially depends on  the following estimate for the volume along
Ricci flow in \cite{Zl}.

\begin{lemma}[Lemma 3.1 in \cite{Zl}]\label{l301} Let $\bar{g}(\bar{t}), \bar{t}\in[0,T),$ be a solution to the Ricci flow equation
(\ref{00*}), i.e. $$\frac{\partial}{\partial
\bar{t}}\bar{g}(\bar{t})=-2Ric_{\bar{t}}$$ on a closed manifold $M$.
If $\lambda_{M}(\bar{g}(\bar{t}))\leq0$,  for all $\bar{t}$, then
there exist constants $c_{1},c_{2}>0$ depending only on
$\bar{g}(0)$, such that,  for all $\bar{t}\geq0$,
$${\rm{Vol}}_{\bar{g}(\bar{t})}(M)\geq c_{1}e^{-c_{2}\bar{t}}.$$
\end{lemma}

 Since the importance of this lemma, we present  the sketch of the
 proof here for reader's  convenience.

\begin{proof}[Sketch of the
 proof]
First, we recall some basics about the  $\mu$ functional introduced
by Perelman \cite{P1}.  Given a closed Riemannian manifold $(M,g)$
and a function $f\in C^{\infty}(M)$ and a constant $\tau>0$, define
$$\mathcal{W}(g,f,\tau)=\int_{M}[\tau(R_{g}+|\nabla
f|^{2})+f-n](4\pi\tau)^{-n/2}e^{-f}dv_{g},$$ and then set
\begin{equation}\label{e2.001}
\mu(g,\tau)=\inf\{\mathcal{W}(g,f,\tau)|\int_{M}(4\pi\tau)^{-n/2}e^{-f}dv_{g}=1\}.
\end{equation}
 By a result of Rothaus \cite{R}, for each
$\tau>0$, there is a smooth minimizer $f$ such that
$\mu(g,\tau)=\mathcal{W}(g,f,\tau)$. In \cite{Zl}, some bounds of
the $\mu$-functional were  obtained.  First, Lemma 2.1  in \cite{Zl}
shows that  there is  a  lower bound for $\tau>\frac{n}{8}$,
\begin{equation}\label{e2.100}
\mu(g,\tau)\geq \lambda_{M}(g)
\tau-\frac{n}{2}\ln(4\pi\tau)-n-\frac{n}{8}(\lambda_{M}(g) -\inf
R_{g})-n\ln C_{s},
\end{equation}
where $C_{s}$ denotes the Sobolev constant for $g$, i.e.
$\|\phi\|_{L^{\frac{2n}{n-2}}(g)}\leq C_{s}\|\phi\|_{H^{1,2}(g)}$
for all $\phi\in C^{\infty}(M)$.
 Second, Corollary  2.3 in \cite{Zl} says that, if $\lambda_{M}(g)\leq0$, then
\begin{equation}\label{e2.101}
\mu(g,\tau)\leq\ln{\rm{Vol}}_{g}(M)-\frac{n}{2}\ln(4\pi\tau)-n+1.
\end{equation}
 In  \cite{P1},
 Perelman proved the monotonicity of $\mu$-functional along the Ricci
flow:

\begin{theorem}[\cite{P1}]\label{p11}
Let $\bar{g}(\bar{t})$ be a solution to the Ricci flow equation
(\ref{00*}) on a closed manifold $M$. Denote
$\tau(\bar{t})=A-\bar{t}$ for some constant $A>0$, then
$\mu(\bar{g}(\bar{t}),\tau(\bar{t}))$ is non-decreasing whenever it
makes sense.
\end{theorem}

Note that   $ \lambda_{M}(\bar{g}(\bar{t}))\leq 0$.  Substituting
$\tau=\frac{n}{8}$ into (\ref{e2.101}), then using Theorem \ref{p11}
and (\ref{e2.100}), we have
\begin{eqnarray}
{\rm{Vol}}_{\bar{g}(\bar{t})}(M)&\geq&\exp(\mu(\bar{g}(\bar{t}),\frac{n}{8})+\frac{n}{2}\ln(\frac{n}{2}\pi)+n-1)\nonumber\\
&\geq&\exp(\mu(\bar{g}(0),\frac{n}{8}+\bar{t})+\frac{n}{2}\ln(\frac{n}{2}\pi)+n-1)\nonumber\\
&\geq&\exp(\lambda_{M}(\bar{g}(0))
\bar{t}-\frac{n}{2}\ln(1+\frac{8}{n}\bar{t})+\frac{n}{8}\inf
R_{0}-n\ln C_{s}(\bar{g}(0))-1)\nonumber\\
&\geq&\exp((\lambda_{M}(\bar{g}(0))-4)\bar{t}+\frac{n}{8}\inf
R_{0}-n\ln C_{s}(\bar{g}(0))-1).
\end{eqnarray}
 We obtain the conclusion.
\end{proof}

Now we can prove:

\begin{lemma}\label{l35} Let $M$ be a closed   $n$-manifold with $\overline{\lambda}_{M}\leq 0 $, and
 $g(t),t\in[0,\infty),$ be   a long time  solution  of (\ref{*}) with scalar curvature  $|R_{t}|<C $ for a constant $C$
 independent of $t$. We
 have
$$\liminf_{t\rightarrow\infty}r(t)=\liminf_{t\rightarrow\infty}\int_M R_{t}dv_{g(t)}\leq 0.$$
\end{lemma}

\begin{proof}  If it is  not true,  there is a
constant $\delta>0$ such that, for $t\gg 1$,  $ r(t)>\delta $. By a
translation on $t$, we assume that $ r(t)>\delta $ for all $t>0$.

Let $\bar{g}(\bar{t})=\sigma (t) g(t), \bar{t} \in [0, T),$ be the
corresponding Ricci flow solution, i.e. $ \frac{\partial}{\partial
\bar{t}}\bar{g}(\bar{t})=-2Ric_{\bar{t}}$ with $\bar{g}(0)=g(0) $,
which implies  $\sigma(t)=\exp (-\frac{2}{n}\int_{0}^{t}r(s)ds) $
and $\bar{t}= \int_{0}^{t}\sigma(s)ds$. Thus
$$T=\int_0^\infty \sigma(s)ds=\int_0^\infty \exp(-\int_0^s\frac{2}{n}r(u)du)ds<
\frac{n}{2\delta}.$$ Now we can compute
\begin{eqnarray}
(T-\bar{t})Vol_{\bar{g}(\bar{t})}(M)^{-\frac{2}{n}}&=& \sigma(t)^{-1}(T-\bar{t})\nonumber\\
&=&\exp(\int_0^t\frac{2}{n}r(u)du)\cdot\int_t^\infty
\exp(-\int_0^s\frac{2}{n}r(u)du)ds\nonumber\\
&\geq &\int_t^\infty\exp(-\int_t^s\frac{2}{n}r(u)du)ds,\nonumber \\
& \geq & \int_t^\infty\exp(-\frac{2}{n}C(s-t))ds,\nonumber \\  &
\geq & \frac{n}{2C},
\end{eqnarray}
 by   $r(t)\leq \sup |R_{t}| \leq C$ for a constant $C>0$ independent of $t$.
Thus  there exists a  $C_1<\infty$  such that
\begin{equation}\label{e2.10} Vol_{\bar{g}(\bar{t})}(M)\leq
C_1(T-\bar{t})^{n/2}, \end{equation} which  implies that
$$\lim_{\bar{t}\rightarrow
T}{\rm{Vol}}_{\bar{g}(\bar{t})}(M)=0.
 $$

However,  since $\overline{\lambda}_{M}\leq 0 $, we have
$\lambda_{M}(\bar{g}(\bar{t}))\leq 0 $, and
$$ {\rm{Vol}}_{\bar{g}(\bar{t})}(M)\geq c_{1}e^{-c_{2}\bar{t}}
\geq c_{1}e^{-c_{2}T}, $$ for two constants $c_{1}>0$ and $c_{2}>0$
by Lemma  \ref{l301}, which is a contradiction.
\end{proof}

Before  proving  Theorem \ref{1.1}, we recall   the evolution
equations for volume forms and scalar curvatures along the
normalized Ricci flow (\ref{*}):
\begin{equation}\label{e2.3}\frac{\partial}{\partial t}dv_{g(t)}=-(R_{t}-r(t))dv_{g(t)}, \ \ \ \ {\rm
and}\end{equation}
\begin{equation}\label{e2.4}\frac{\partial}{\partial t}R_{t}= \triangle_{t}
R_{t}+2|Ric_{t}\textordmasculine|^{2}+\frac{2}{n}R_{t}(R_{t}-r(t)),
\end{equation} where
$Ric_{t}\textordmasculine=Ric_{t}-\frac{R_{t}}{n}g(t)$ denotes the
Einstein tensor (c.f. \cite{H2} ).

\begin{proof}[Proof of Theorem \ref{1.1}] Note that $$\breve{R}_{t}=\inf_{M}R_{t}\leq \lambda_{M}(g(t))\leq
\overline{\lambda}_{M}\leq 0.$$
 If $\breve{R}_{t}\leq -c <0 $ for a constant $c>0$, we obtain the conclusion
 by Lemma 2.7  and
     Lemma
 3.1 in
\cite{FZZ1}.

From the maximal principal, $$\frac{\partial}{\partial
t}\breve{R}_{t}\geq \frac{1}{2}\breve{R}_{t}(\breve{R}_{t}-r(t))\geq
0, $$ and thus $\breve{R}_{t}$ is non-decreasing.   So only case
left to prove is that
$$\lim_{t\longrightarrow \infty}\breve{R}_{t}=0.  $$
 By Lemma \ref{l35}, and $\breve{R}_{t}\leq r(t)$,
$$\liminf_{t\rightarrow\infty}r(t)=0, $$ which implies that  $\overline{\lambda}_{M}=0 $.

 First, we assume that there is
a sequence $t_{k}'\longrightarrow \infty$ such that
$r(t_{k}')>\epsilon$ for a constant $\epsilon >0$ independent of
$k$.  Since $\liminf_{t\rightarrow\infty}r(t)=0$, there are
$t_{k}\in (t_{k}',t_{k+1}')$ such that, for $k\gg 1$,
 $$r(t_{k})=\inf_{(t_{k}',t_{k+1}')}r(t)\longrightarrow 0, \ \ \ \ \frac{dr}{dt}
 (t_{k})=0.$$

Now, we  assume $ \lim_{t\rightarrow\infty}r(t)=0$.  If
$|\frac{dr}{dt}(t)|>\delta >0 $ for a constant $ \delta$ independent
of $t$  when $t\gg 1$,  $ |r(t)-r(0)|> \delta t$ which is a
contradiction. Thus there is a sequence  $t_{k}\longrightarrow
\infty$ such that
 $$\lim_{t_{k} \longrightarrow \infty} r(t_{k})=0, \ \ \ \ \lim_{t_{k} \longrightarrow \infty}\frac{dr}{dt} (t_{k})=0.
$$

In both cases, we have  \begin{eqnarray}
0&=&\lim_{k\rightarrow\infty}\frac{dr}{dt}(t_k)\nonumber\\
&=&\lim_{k\rightarrow\infty}\int_M(2|Ric^o_{t_k}|^2-\frac{1}{2}R_{t_k}(R_{t_k}-r(t_k)))dv_{g(t_k)}\nonumber\\
&\geq&\lim_{k\rightarrow\infty}\int_M2|Ric^o_{t_k}|^2dv_{g(t_k)}-\lim_{k\rightarrow\infty}C\int_M|R_{t_k}-r(t_k)|dv_{g(t_k)}\nonumber\\
&\geq&\lim_{k\rightarrow\infty}\int_M2|Ric^o_{t_k}|^2dv_{g(t_k)}-\lim_{k\rightarrow\infty}C\int_M(R_{t_k}+r(t_k)-2\breve{R}_{t_k})
dv_{g(t_k)}\nonumber\\
&=&\lim_{k\rightarrow\infty}\int_M2|Ric^o_{t_k}|^2dv_{g(t_k)}-\lim_{k\rightarrow\infty}2C(r(t_k)-\breve{R}_{t_k})\nonumber\\
&=&\lim_{k\rightarrow\infty}\int_M2|Ric^o_{t_k}|^2dv_{g(t_k)},\nonumber
\end{eqnarray} by (\ref{e2.3}), (\ref{e2.4}), and   the assumption $| R_{t}|< C$ for a constant $C$
independent of $t$.

The Chern-Gauss-Bonnet formula and the Hirzebruch signature theorem
(cf. \cite{Be}) show  that,  for any metric $g$ on $M$,
$$\chi(M)=
\frac{1}{8\pi^{2}}\int_{M}(\frac{R_{g}^{2}}{24}+|W^{+}_{g}|^{2}+|W^{-}_{g}|^{2}-\frac{1}{2}
|Ric_{g}\textordmasculine|^{2})dv_{g}, \ \ \ \ { \rm and}$$
$$\tau(M)=\frac{1}{12\pi^{2}}\int_{M}(|W^{+}_{g}|^{2}-|W^{-}_{g}|^{2})dv_{g},$$
where $W^{+}_{g}$ and $W^{-}_{g}$ are the self-dual and
anti-self-dual Weyl tensors respectively.  Thus
\begin{eqnarray*} 2\chi(M)-3|\tau(M)|& \geq &
\liminf\limits_{k\longrightarrow \infty }
\frac{1}{4\pi^{2}}\int_{M}(\frac{1}{24}
R_{t_{k}}^{2}-\frac{1}{2}|Ric_{t_{k}}\textordmasculine|^{2})dv_{g(t_{k})}
\\ & = & \liminf\limits_{k\longrightarrow \infty }
\frac{1}{4\pi^{2}}\int_{M}\frac{1}{24}R_{t_{k}}^{2} dv_{g(t_{k})}
\ge 0.\end{eqnarray*} Since $\overline{\lambda}_{M}=0$ in this case,
we obtain the conclusion.
\end{proof}

\end{document}